\documentclass[14pt]{article}
\usepackage{ amssymb }
\usepackage{lipsum}

\usepackage{amsfonts}
\usepackage{amsthm}

\usepackage{epsfig}
\usepackage{graphicx}  
  
\usepackage{pdfsync}

\usepackage{caption}
\usepackage{subcaption}
 
\usepackage{amsmath,amsfonts,amssymb,amsthm,epsfig,epstopdf,titling,url,array}

\usepackage{hyperref}

 \usepackage{tikz}

\theoremstyle{plain}
\newtheorem{thm}{Theorem}[section]
\newtheorem{lem}[thm]{Lemma}
\newtheorem{prop}[thm]{Proposition}
\newtheorem*{cor}{Corollary}

\theoremstyle{definition}
\newtheorem{defn}{Definition}[section]

\newtheorem{exmp}{Example}[section]

\theoremstyle{remark}
\newtheorem*{rem}{Remark}

\theoremstyle{remark}

\usepackage{theoremref}

\usepackage{thmtools}

\usepackage{graphicx}
\usepackage{xcolor}
\definecolor{darkblue}{rgb}{0,0,0.5} 
\usepackage{transparent}

\DeclareMathOperator{\diam}{diam}

\newtheorem{example}[thm]{Example}

\newtheorem{theorem}{Theorem}[section]

\begin{document}

\title{\textbf{Slim Polygons, Hyperbolicity, and a new Quasi-Isometric Invariant}}
\author{Panagiotis Tselekidis}

\maketitle


\begin{abstract}
We introduce the slimness function, $\mathrm{Slim}_X(n)$,  of a geodesic metric space, which
records the maximal slimness of its geodesic polygons as a function of
the number of sides. We study the large-scale behaviour of this function and its relation to hyperbolicity, growth function, and quasi-isometries.

We first establish that the slimness
function is a natural quasi-isometric invariant in the hyperbolic
setting. We then obtain quantitative estimates for geodesic polygons
in hyperbolic spaces and bounded-valency graphs. In particular, we
show that sufficiently slow growth of the slimness function forces
tree-like geometry. We also prove a decomposition result for trees of hyperbolic spaces.

We show that for a finitely generated hyperbolic group
$G$ there are only two possibilities: either $G$ is virtually free and
$\operatorname{Slim}_G$ is bounded, or $\operatorname{Slim}_G(n)\simeq\log n$.

Finally, we investigate the interaction between polygonal slimness and
the growth of connected bounded-valency graphs. Combining the resulting
estimates with large-scale dimension theory, we obtain criteria
implying that an unbounded graph has asymptotic dimension one. 
\end{abstract}

\tableofcontents

\section{Introduction} The geometry of metric spaces, particularly in the realm of geometric group theory, is deeply intertwined with large-scale invariants that capture coarse geometric structures. Classical notions such as Gromov hyperbolicity rely heavily on the concept of $\delta$-slim triangles (and polygons). In this paper, we systematically investigate the quantitative aspects of slimness for polygonal structures across various geodesic metric spaces.

Our primary object of study is the \emph{slimness function}, which measures the ``optimal'' constant for which every geodesic $n$-gon in a space $X$ is $\alpha$-slim. We prove that $\mathrm{slim}_X$ behaves as a quasi-isometric invariant when the underlying space $X$ is hyperbolic. Moreover, we analyze how this function behaves under finite graph decompositions of groups, showing that $\mathrm{slim}_G$ is controlled by the supremum over the vertex groups. A central motivation of our work is the interaction between slimness, growth functions, and asymptotic dimension ($\mathrm{asdim}$).



It is natural to ask how much information about the
large-scale geometry of a space is contained in the corresponding
slimness properties of geodesic polygons with an arbitrary number of
sides.

For a geodesic metric space $X$ and an integer $n\geq 3$, let
$\operatorname{Slim}_X(n)$ denote the supremum of the slimness
constants of geodesic $n$-gons in $X$.
Let $P$ be a geodesic
$n$-gon in $X$. We define $\operatorname{Slim_X}(P)$ to be the infimum
of all $\delta\geq 0$ such that $P$ is $\delta$-slim. 
 Thus, the function

$$
n\longmapsto \operatorname{Slim}_X(n)
$$

records quantitatively how the geometry of geodesic polygons
deteriorates as the number of sides increases. In a Gromov hyperbolic
space (see \cite{Gr}), geodesic triangles are uniformly slim, and consequently every
geodesic $n$-gon is $(n-2)\delta$-slim for a suitable hyperbolicity
constant $\delta$. Hence the function $\operatorname{Slim}_X(n)$ is at
most linear in $n$ in every hyperbolic space. On the other hand, this
linear estimate is obtained by decomposing an $n$-gon into geodesic
triangles and does not necessarily reflect the actual behaviour of
polygonal slimness.

The purpose of this paper is to study the function
$\operatorname{Slim}_X$ as a large-scale geometric invariant. Rather
than considering only whether geodesic polygons are uniformly slim,
we investigate the asymptotic rate at which their slimness may grow
with the number of sides. 

A first difficulty is that the precise function
$\operatorname{Slim}_X(n)$ is not expected to be invariant under
quasi-isometries. Quasi-isometries do not preserve geodesics exactly. Consequently, polygonal slimness can change under a
quasi-isometry by multiplicative and additive errors. This suggests
that the appropriate object is not the function itself, but rather
its asymptotic growth type.

For two functions $f,g\colon\mathbb N\to[0,\infty)$, we write
$f\preceq g$ if there exist constant $C>0$ such that

$$
f(n)\leq Cg(Cn+C)+C
$$

for all sufficiently large $n$. We then may say that $f\preceq g$ with constant $C $. We write $f\simeq g$ if
$f\preceq g$ and $g\preceq f$.
We write that $f\simeq g$ with constant $C $, if both $f\preceq g$ and $g\preceq f$ are satisfied with constant $C $. The first main result of this paper
shows that the asymptotic growth type of the polygonal slimness
function is preserved by quasi-isometries on hyperbolic spaces.

\begin{theorem}
Let $X$ and $Y$ be geodesic hyperbolic metric spaces. If $X$ and $Y$ are
quasi-isometric, then

$$
\operatorname{Slim}_X\simeq \operatorname{Slim}_Y.
$$

\end{theorem}

We give the following examples:
\begin{example}
If $X$ is a hyperbolic space, then $\operatorname{Slim}_X \preceq n$ (\cite{Gr}).
\end{example}

\begin{example}
Let $X$ be an unbounded vertex-transitive connected hyperbolic graph of valency at most $k$. Then $\operatorname{Slim}_X \preceq \sqrt{ \log_k (n)}$ if and only if $X$ is quasi-tree (\cite{PT Polygon}).
\end{example}
\begin{example}
Let $G$ be an unbounded connected $\delta$- hyperbolic graph of valency at most $k$. Then $\operatorname{Slim}_X \leq \frac{1}{8}\log_k(n\delta)\text{-slim}$ if and only if $X$ is quasi-tree (\cite{PT Polygon}).

\end{example}

\begin{example}
In Euclidean spaces, if $P$ is a
geodesic $n$-gon in a Euclidean space and
$l_0\leq \cdots\leq l_n$ are the side lengths of $P$, then $\operatorname{Slim}_X(P) \leq (n-2)l_{n-2}$.
\end{example}

\begin{example}

Let $X$ be a geodesic graph with valency at most $k$. Suppose that for every geodesic $n$-gon, with side lengths $l_0 \leq\cdots\leq l_{n}, $ we have $ \operatorname{Slim}_X(P) \leq \frac{1}{8}\log_k\left((n-2)\frac{l_{n-2}}{17}\right).$ Then
$\operatorname{asdim}_{AN}X=1$ (\cite{PT Polygon}).

\end{example}

It is natural to ask the following

\begin{flushleft}
\textbf{Question 1:}\label{question:main}
To what extent does the asymptotic growth of
$\operatorname{Slim}_X(n)$ determine the large-scale geometry of
$X$?
\end{flushleft}

 We next consider the
behaviour of this invariant under trees of spaces decompositions. In
particular, we study the slimness of hyperbolic spaces which split as trees of hyperbolic spaces. We give the following
\begin{theorem}\label{thm.trees of spaces}
Let $\delta\geq 0$ and $M\geq 0$. Assume that $X$ is a tree of
$\delta$-hyperbolic spaces with finite, $M$-bounded edge spaces. Assume
moreover that there are geodesic metric spaces $Y_1,\dots,Y_m$ and constants
$\lambda\geq1$, $\varepsilon\geq0$ such that every vertex space $X_v$ is
$(\lambda,\varepsilon)$-quasi-isometric to some $Y_j$, and every $Y_j$ is
$(\lambda,\varepsilon)$-quasi-isometric to some vertex space $X_v$. Then
$$
\operatorname{Slim}_X
\simeq
\max_{1\leq j\leq m}\operatorname{Slim}_{Y_j}.
$$
\end{theorem}

Observe that each $Y_j$ is hyperbolic, being a geodesic metric space
quasi-isometric to a $\delta$-hyperbolic space.

Geodesic polygons with many sides in hyperbolic spaces are a classical
tool, they were used extensively by Ol'shanskii in his work on periodic
quotients and small cancellation theory over hyperbolic groups
\cite{Ol91,Ol93}. The point of view of \cite{PT Polygon} and this paper is
different: one asks how fast $\operatorname{Slim}_X(n)$ may grow and what
this growth says about the large-scale geometry of $X$. In \cite{PT Polygon}
it is shown that a sufficiently small logarithmic bound on the slimness of
polygons in a graph of bounded valency forces the graph to be a
quasi-tree, and that a vertex-transitive hyperbolic graph of bounded
valency with $\operatorname{Slim}_X\preceq\sqrt{\log_k(n)}$ is a
quasi-tree. The following theorem determines $\operatorname{Slim}_G$ for every hyperbolic
group $G$.

\begin{theorem}\label{dich:thmA}
Let $G$ be a finitely generated hyperbolic group.
\begin{enumerate}
\item If $G$ is virtually free, then $\operatorname{Slim}_G$ is bounded.
\item If $G$ is not virtually free, then $\operatorname{Slim}_G(n)\simeq\log n$.
\end{enumerate}
\end{theorem}

The upper bound $\operatorname{Slim}_G(n)\preceq\log n$ holds in every
hyperbolic space, we include a short direct proof with an
explicit constant (Proposition~\ref{dich:upper}). The lower bound combines
a  construction in the hyperbolic plane
(Proposition~\ref{dich:lower}) with the theorem of Bonk and Kleiner
\cite{BK} that every hyperbolic group which is not virtually free contains
a quasi-isometrically embedded copy of $\mathbb{H}^2$.

Since a group with $\operatorname{Slim}_G(3)<\infty$ is hyperbolic,
Theorem~\ref{dich:thmA} gives the following characterisation of virtually
free groups.

\begin{cor}[Characterisation of virtually free groups]
Let $G$ be a finitely generated group with a finite generating set. The
following are equivalent:
\begin{enumerate}
\item $G$ is virtually free;
\item $\operatorname{Slim}_G$ is bounded;
\item $\operatorname{Slim}_G(n)<\infty$ for all $n$ and
$\operatorname{Slim}_G(n)=o(\log n)$.
\end{enumerate}
\end{cor}

\begin{cor}[Surface groups]
Let $G$ be the fundamental group of a closed orientable surface of genus
$g\geq 2$, or more generally a group acting properly and cocompactly by
isometries on $\mathbb{H}^m$, $m\geq 2$. Then
$\operatorname{Slim}_G(n)\simeq\log n$.
\end{cor}

This corollary follows directly from Theorem~\ref{dich:thmA}, and also
from Proposition~\ref{dich:Hm} together with the quasi-isometry invariance
of the slimness function (Theorem~\ref{thm3.1}).

Thus, among hyperbolic groups, the growth type of the slimness function
takes exactly two values, and the logarithmic upper bound valid in every
hyperbolic space is attained by every hyperbolic group which is not
virtually free. We do not know whether a similar dichotomy holds for
hyperbolic graphs of bounded valency; see Section~\ref{dich:questions}.

\medskip

This theorem provides a useful
mechanism for transferring information about slimness from the
building blocks of a group to the group itself.

We next turn to the relation between slimness and growth.
For a connected graph $X$ of bounded valency, let

$$
\beta_X(r)=\sup_{x\in X}|B_X(x,r)|
$$

denote its growth function. 

The following result establishes a direct connection between the growth function of a graph, the slimness of its
geodesic polygons, and its asymptotic dimension. More precisely, if the growth function of $X$ is
asymptotically described by an invertible function $f$, then a
sufficiently strong upper bound on the slimness of an $n$-gon, expressed
in terms of $f^{-1}$ and the second largest side length, forces the
space to have asymptotic dimension one. Thus, the
growth--slimness relation provides a quantitative criterion for
one-dimensional large-scale geometry.

\begin{theorem}
Let $X$ be an unbounded vertex-transitive connected graph. Suppose that there exist an invertible, non-decreasing,  $f$ and a constant $C>0$ such that $\beta_X(n) \simeq f(n)$ with constant $C$. We further suppose that $\beta_X(n) \succ n$. We assume that for every geodesic polygon with side lengths $l_0 \leq\cdots\leq l_{n}, $ we have 

$$
\operatorname{Slim}_X(P)
\leq \phi
f^{-1}(n  l_{n-2}),
$$

where $\phi \leq \frac{1}{32C}$. Then

$$
\operatorname{asdim}X=1.
$$

\end{theorem}

The significance of this theorem is that it does not assume
hyperbolicity. Instead, it compares two independent quantitative
features of the graph: the rate at which metric balls grow and the
rate at which the slimness of geodesic polygons is allowed to
increase. The conclusion shows that sufficiently slow polygonal
thickening, relative to the growth of the graph, forces the space to have one-dimensional asymptotic geometry. As a corollary of this theorem we have 

\begin{cor}
Let $X$ be an unbounded vertex-transitive connected graph. Suppose that there exist $r>1$ such that $\beta_X(n) \simeq n^r$ with constant $C$. We assume that for every geodesic polygon with side lengths $l_0 \leq\cdots\leq l_{n}, $ we have 

$$
\operatorname{Slim}_X(P)
\leq 
 \frac{1}{32C}\sqrt[r]{nl_{n-2} }.
$$
 Then

$$
\operatorname{asdim}_{}X=1.
$$
\end{cor}

The paper is organized as follows. In Section~\ref{prelims}, we
introduce the basic definitions of slim geodesic polygons. In
Section~\ref{section3}, we define the
slimness function, and we prove its quasi-isometry invariance in the hyperbolic setting.  In Section~\ref{section4},
we investigate the interaction between polygonal slimness, the  growth function, and the asymptotic dimension. 
In Section~\ref{section5},
we study polygonal slimness for groups admitting finite
graph-of-groups decompositions and establish the corresponding
upper bounds in terms of the vertex groups. In Section~\ref{section 6} we present several open questions.

\noindent\textbf{Acknowledgements.} The author thanks Goulnara
Arzhantseva for bringing the work of Ol'shanskii \cite{Ol91,Ol93} to his
attention.

\section{Preliminaries}\label{prelims}

If $(X,d)$ is a metric space and $A\subseteq X$, we write
$N_r(A)=\{x\in X:d(x,A)\le r\}$, for the closed $r$-neighbourhood of $A$.\\
Let \(X\) be a metric space and let \(R>0\). An \emph{\(R\)-path} in \(X\) is a finite sequence of points $ x_0,x_1,\ldots,x_m\in X $ such that $ d(x_{i-1},x_i)\le R$ for every $i=1,\ldots,m. $ We say that a subset $Z$ of $X$ is $R-$connected id every two points in $Z$ are connected by an $r$-path. Let $C \subseteq X$, we say that $C$ is an \emph{$R-$component} of $X$ if it is $R$-connected, and for every other $R$-connected subset $Z$, of $X$ such that $C \subseteq Z$, we have that $C=Z$.

The \emph{length} of the \(R\)-path is defined by

$$
\ell(x_0,\ldots,x_m)
=\sum_{i=1}^{m}d(x_{i-1},x_i).
$$

The \emph{\(R\)-length} of the path is the number of \(R\)-edges, namely

$$
\ell_R(x_0,\ldots,x_m)=m.
$$

We say that the \(R\)-path joins \(x_0\) to \(x_m\).

\subsection{Geodesic metric spaces and graphs}

A metric space $(X,d)$ is called \emph{geodesic} if for every $x,y\in X$ there exists an isometric embedding $ \gamma\colon [0,d(x,y)]\longrightarrow X$, such that $\gamma(0)=x$ and $\gamma(d(x,y))=y$.
The image of such an embedding is called a \emph{geodesic} from $x$ to $y$ and is denoted by $[x,y]$. A connected graph will always be regarded as a metric space by giving every edge length one and considering the induced path metric. In particular, graphs are geodesic metric spaces.

For a graph $G$, we denote its valency by

$$
\operatorname{val}(G)=\sup_{v\in V(G)}\deg(v).
$$

We say that $G$ has \emph{bounded valency} if $\operatorname{val}(G)<\infty$.

\subsection{Morse lemma}

We need the following very useful lemma for hyperbolic spaces.
\begin{lem}[Morse Lemma]\label{Morse}
Let $X$ be a geodesic $\delta$-hyperbolic metric space. For every
$\lambda\geq 1$ and $\varepsilon\geq 0$, there exists a constant
$R=R(\delta,\lambda,\varepsilon)\geq 0$ such that every
$(\lambda,\varepsilon)$-quasi-geodesic segment in $X$ is contained in
the $R$-neighbourhood of any geodesic segment joining its endpoints.
Conversely, every such geodesic segment is contained in the
$R$-neighbourhood of the quasi-geodesic.
\end{lem}

\subsection{Geodesic polygons and slimness}

Let $X$ be a geodesic metric space. A \emph{geodesic $n$-gon} is a cyclic sequence

$$
P=[x_1,x_2]\cup [x_2,x_3]\cup\cdots\cup[x_n,x_1],
$$

where each $[x_i,x_{i+1}]$ is a geodesic.

For $\delta\ge0$, a geodesic $n$-gon $P$ is called \emph{$\delta$-slim} if every side of $P$ is contained in the $\delta$-neighbourhood of the union of the remaining sides. Thus, for every $i$, $[x_i,x_{i+1}]
\subseteq
N_\delta\left(
\bigcup_{j\ne i}[x_j,x_{j+1}]
\right).
$

A geodesic metric space is called \emph{$\delta$-hyperbolic} if every geodesic triangle is $\delta$-slim. It is called \emph{hyperbolic} if it is $\delta$-hyperbolic for some $\delta\ge0$.

Notice that hyperbolicity only imposes a uniform bound on geodesic triangles. The corresponding bound for general geodesic polygons need not be uniform in the number of sides. Indeed, a standard consequence of hyperbolicity is that an $n$-gon is $\varepsilon_n$-slim, where $\varepsilon_n$ depends on $n$ and $\delta$, where $
\varepsilon_n\le (n-2)\delta.
$



\subsection{Asymptotic Assouad--Nagata dimension}

We recall Gromov's asymptotic dimension. The \emph{asymptotic dimension} $\operatorname{asdim} X$ of a metric space $X$ is defined as follows. We say that $\operatorname{asdim} X\le n$
if and only if, for every $R>0$, there exists a uniformly bounded cover $\mathcal U$ of $X$ whose $R$-multiplicity is at most $n+1$; that is, every ball of radius $R$ in $X$ intersects at most $n+1$ elements of $\mathcal U$.

The \emph{asymptotic Assouad--Nagata dimension} $\operatorname{asdim}_{AN} X$ of a metric space $X$ is defined similarly, with the additional requirement of linear control on the diameters of the elements of the cover. More precisely,

$$
\operatorname{asdim}_{AN} X\le n
$$

if and only if there exist  constants $C,R_0 >0$ such that, for every $R>R_0$, there exists a cover $\mathcal U$ of $X$ with $R$-multiplicity at most $n+1$ and

$$
\operatorname{diam}(U)\le CR
$$

for every $U\in\mathcal U$. It is easy to see that $ \operatorname{asdim} X \leq \operatorname{asdim}_{AN} X$. The terminology and basic theory of asymptotic Assouad--Nagata dimension can be found in \cite{DS07}.

\section{Slimness function}\label{section3}
For a geodesic metric space $X$ and an integer $n\geq 3$, let $P$ be a geodesic
$n$-gon in $X$. We define $\operatorname{Slim_X}(P)$ to be the infimum
of all $\delta\geq 0$ such that $P$ is $\delta$-slim.
Then we set $$\operatorname{Slim}_X(n)= \sup \{ \operatorname{Slim_X}(P) \mid P \text{ is a geodesic n-gon of X} \}.$$
 
Since degenerate sides are allowed, $\operatorname{Slim}_X(n)$ is
non-decreasing in $n$: adding a side of length zero at a vertex turns a
geodesic $n$-gon into a geodesic $(n+1)$-gon with the same slimness. The \emph{slimness function}

$$
n\longmapsto \operatorname{Slim}_X(n)
$$

records quantitatively how the geometry of geodesic polygons
deteriorates as the number of sides increases. In a Gromov hyperbolic
space (see \cite{Gr}), geodesic triangles are uniformly slim, and consequently every
geodesic $n$-gon is $(n-2)\delta$-slim for a suitable hyperbolicity
constant $\delta$. Hence the function $\operatorname{Slim}_X(n)$ is at
most linear in $n$ in every hyperbolic space. On the other hand, this
linear estimate is obtained by decomposing an $n$-gon into geodesic
triangles and does not necessarily reflect the actual behaviour of
polygonal slimness.

For two functions $f,g\colon\mathbb N\to[0,\infty)$, we write
$f\preceq g$ if there exist constant $C>0$ such that

$$
f(n)\leq Cg(Cn+C)+C
$$

for all sufficiently large $n$. We then may say that $f\preceq g$ \emph{is satisfied with constant $C $}. We write $f\simeq g$ if
$f\preceq g$ and $g\preceq f$.
We write that $f\simeq g$ with constant $C $, if both $f\preceq g$  and $g\preceq f$ are satisfied with constant $C $.

\begin{prop}\label{prop3.1}
Let $X$ and $Y$ be geodesic hyperbolic metric spaces. Assume that $f: X \to Y$ is a $(\lambda,\varepsilon)-$ quasi-isometric embedding. Then $$
\operatorname{Slim}_X \preceq \operatorname{Slim}_Y.
$$
\end{prop}

\begin{proof}
Let $P$ be a geodesic n-gon of $X$. Then $ P=[x_1,x_2]\cup [x_2,x_3]\cup\cdots\cup[x_n,x_1],$ where each $r_i=[x_i,x_{i+1}]$ is a geodesic. We observe that the image $f(P)$ may no longer be a geodesic n-gon in $Y$ (since $f(r_i)$'s are $(\lambda,\varepsilon)-$ quasi-geodesics in $Y$).\\
We now consider $\gamma_i=[f(x_i),f(x_{i+1})]$ geodesics in $Y$. Those geodesics form a geodesic $n$-gon, say $Q$, in $Y$. 
Let $x\in r_i$, for some $i$.

By lemma \ref{Morse} there exists a constant
$R=R(\delta,\lambda,\varepsilon)\geq 0$ such that every
$(\lambda,\varepsilon)$-quasi-geodesic $f(r_i)$ is contained in
the $R$-neighbourhood of $\gamma_i$.
So, there exists $y  \in \gamma_i$, such that $d_Y(f(x),y) \leq R$.\\
Since $Q$, is $\operatorname{Slim}_Y(n)-$slim, there exists $y_1 \in \gamma_{j} $ ($i \neq j$) such that $d_Y(y,y_1) \leq \operatorname{Slim}_Y(n)$. By Morse lemma again, there exists, $x_1 \in r_{j} $ such that $d_Y(y_1,f(x_1))\leq R$.
 Thus $d_Y(f(x),f(x_1)) \leq \operatorname{Slim}_Y(n) + 2R(\delta,\lambda,\varepsilon)$. We conclude that $$d_X(x,x_1) \leq \lambda (\operatorname{Slim}_Y(n) + 2R(\delta,\lambda,\varepsilon) +\varepsilon).$$

Since $x_1\in r_j$ with $j\neq i$, and $x\in r_i$ was arbitrary, it follows that
$\operatorname{Slim}_X(n)\leq\lambda\operatorname{Slim}_Y(n)+\lambda(2R+\varepsilon)$
for all $n$. Hence $\operatorname{Slim}_X \preceq \operatorname{Slim}_Y$.

\end{proof}

As a corollary we have that the slimness
function is preserved by quasi-isometries on hyperbolic spaces.

\begin{theorem}\label{thm3.1}
Let $X$ and $Y$ be geodesic hyperbolic metric spaces. If $X$ and $Y$ are
quasi-isometric, then

$$
\operatorname{Slim}_X\simeq \operatorname{Slim}_Y.
$$

\end{theorem}

We give the following examples:
\begin{example}
If $X$ is a hyperbolic space, then $\operatorname{Slim}_X \preceq n$ (\cite{Gr}).
\end{example}

\begin{example}
Let $X$ be an unbounded connected hyperbolic graph of valency at most $k$. Then $\operatorname{Slim}_X \preceq \sqrt{ \log_k (n)}$ if and only if $X$ is quasi-tree (\cite{PT Polygon}).
\end{example}
\begin{example}
Let $G$ be an unbounded connected $\delta$- hyperbolic graph of valency at most $k$. Then $\operatorname{Slim}_X \leq \frac{1}{8}\log_k(n\delta)\text{-slim}$ if and only if $X$ is quasi-tree (\cite{PT Polygon}).

\end{example}

\begin{example}
In Euclidean spaces, if $P$ is a
geodesic $n$-gon in a Euclidean space and
$l_0\leq \cdots\leq l_n$ are the side lengths of $P$, then $\operatorname{Slim}_X(P) \leq (n-2)l_{n-2}$.
\end{example}

\begin{example}

Let $X$ be a geodesic graph with valency at most $k$. Suppose that for every geodesic $n$-gon, with side lengths $l_0 \leq\cdots\leq l_{n}, $ we have $ \operatorname{Slim}_X(P) \leq \frac{1}{8}\log_k\left((n-2)\frac{l_{n-2}}{17}\right).$ Then
$\operatorname{asdim}_{AN}X=1$ (\cite{PT Polygon}).

\end{example}

\subsection{Euclidean Spaces}
We will show that, if $P$ is a
geodesic $n$-gon in a Euclidean space $X$ and
$l_1\leq\cdots\leq l_n$ are the side lengths of $P$, then  

\begin{center}
$ \frac{ (n-2)l_{n-2}}{2} \leq \operatorname{Slim_X} P \leq (n-2)l_{n-2}$
\end{center}

 It suffices to show that for the plane. For the upper bound, it is not hard to prove it for triangles. Let $P$ be a triangle with sides $a \leq b \leq c$. Then for every point $x$ of $b$ we consider a parallel line to $a$, which intersects $c$ at a point   $y$. Then by elementary Euclidean geometry we know that $d(x,y)\leq a$.\\ 
By using induction on the number of vertices of the n-gons, and the fact that we can separate the two largest sides of an n-gon P by a geodesic we can obtain the wanted upper bound.

For the lower bound, for any $n>2$, we consider a triangle $P$, with sides $a=b <c$, with a right angle between the sides $a$ and $b$. We decompose $a$ into $n-2$ segments of equal length, say $\lambda$. Then we have an n-gon, say $Q$, with sides $l_1,\dots,l_{ n-2}=\lambda<b<c$. If $x$ is the mid point of $c$, the distance of $x$ from the other sides is $\frac{b}{2}=\frac{(n-2)\lambda}{2} = \frac{(n-2) l_{n-2}}{2}$.

\section{Slimness and Growth functions}\label{section4}
The next theorem is the main applications of the growth-slimness
relation developed in this paper. It shows that, for bounded-valency
graphs, a sufficiently strong control of the slimness of geodesic
polygons in terms of the inverse of the growth function forces
one-dimensional large-scale geometry. In particular, the theorem gives
a quantitative criterion for
$\operatorname{asdim}X=1$, linking the combinatorial growth of the
graph directly to the geometry of its geodesic polygons.

\begin{theorem}
Let $X$ be an unbounded vertex-transitive connected graph. Suppose that there exist an invertible, non-decreasing,  multiplicative function $f$ and a constant $C>0$ such that $\beta_X(n) \simeq f(n)$ with constant $C$. We further suppose that there is $k>0$ so that $f(x)\geq kx$, for all $x\geq x_0$. We assume that for every geodesic polygon with side lengths $l_0 \leq\cdots\leq l_{n}, $ we have 

$$
\operatorname{Slim}_X(P)
\leq \phi 
f^{-1}(n l_{n-2}),
$$

where $\phi=f^{-1}(\frac{k}{64C^3} \frac{1}{2 f(32C)})$. Then

$$
\operatorname{asdim}X=1.
$$
\end{theorem}
\begin{proof}
We may assume that $C\geq 10$. Then we have $
f(n)\leq C\beta_X(Cn+C)+C$ and  $\beta_X(n)\leq Cf(Cn+C)+C$, for all $n\geq n_0.$\\
Let $R> \max \{ \frac{8C^2}{k} ,2C(x_0 +1) ,4C, 2n_0 \}$. Let $e\in X$ be the base point. We decompose $X$ using annuli $A_i=B_{3R(i+1)}(e)\setminus B_{3Ri}(e)$, where $i \in \mathbb{N}$. We will construct suitable covers, say $\mathcal{U}_i$, for  annuli $A_i$.
For the annulus $A_i$, we consider its R-components.

\begin{flushleft}
\textbf{Claim:} If $C$ is an R-component of $A_i$ ($i\geq i_0$), then $\operatorname{diam}C \leq mR.$ Where $m=30$.
\end{flushleft}
 Set $i_0=10[2R+4\phi f^{-1}(4R \beta_X((m+2)R))]$.
Suppose for a contradiction that $\operatorname{diam}C >mR$. Then there are $x,y \in C$ such that $$mR+1\leq d(x,y) \leq (m+1)R .$$

Let $x=x_0, \dots, x_{n-2}=y$ be an R-path from $x$ to $y$ of minimal R-length, we denote this path by $\gamma$. There exists a number $ 0<\epsilon<1$ such that 
\begin{equation}\label{eq1}
d(x,y)=mR+ \epsilon R.
\end{equation}  Let $r_x$, $r_y$ be geodesics from $e$ to $x$, and from $e$ to $y$, respectively. Observe that the geodesics $r_x, [x_0,x_1], \dots, [x_{n-1},x_n], r_y $ form a geodesic $n-$gon.

 From the hypothesis, this polygon is $\phi f^{-1}(n l_{n-2})\text{-slim}.$, so $x'=r_x(3Ri-(\frac{3R}{2}+ \phi
f^{-1}(n l_{n-2}))$ is $\phi
f^{-1}(n l_{n-2})$-close to one of the other geodesics.

From the way we choose $x'$, we can see that $\operatorname{dist}(x',[x_j,x_{j+1}])>\phi
f^{-1}(n l_{n-2})$. Thus, $\operatorname{dist}(x',r_y) \leq \phi
f^{-1}(n l_{n-2})$.
In particular, we can deduce that $d(x',y')\leq 2\phi
f^{-1}(n l_{n-2})$, where $y'=r_y(3Ri-(\frac{3R}{2}+\phi
f^{-1}(n l_{n-2})))$. 

We have $d(x,y)\leq 9R +  4 \phi
f^{-1}(n l_{n-2})$, we deduce that

\begin{equation}\label{eq13}
(m-9)R\leq  4 \phi
f^{-1}(n l_{n-2})
\end{equation}

We set $\overline{\gamma}$ to be the path formed by a union of the geodesics  $ [x_0,x_1], \dots, [x_{n-1},x_n] $.\\ 
\textbf{Observation 1}\label{ob1}
Since $\gamma$ has minimal R-length we observe that the set $\{ x_{2k} \mid k \in \mathbb{N}\}$ is $\frac{R}{2}$-disjoint.\\

We distinguish two cases:
\begin{flushleft}
\textbf{Case 1:} $\gamma \subseteq B_{d(x,y)}(x)$.\\
Then $\overline{\gamma} \subseteq B_{d(x,y)+\frac{R}{2}}(x)$. Thus by observation 1 there are at least $\frac{n}{2}-1$ many disjoint $\frac{R}{2}-$balls in $B_{d(x,y)+\frac{R}{2}}(x).$ Thus, by vertex-transitivity there are at least $(\frac{n}{2}-1) \beta_X (\frac{R}{2})$ many elements in $B_{d(x,y)+\frac{R}{2}}(x).$ So,  

\begin{equation}\label{eq15}
(\frac{n}{2}-1)\beta_X (\frac{R}{2}) \leq \beta_X ((m+2)R).
\end{equation}

We note that $n\geq 8$, since otherwise $d(x,y)\leq8R$, but  $d(x,y) \geq mR>20R$.

Thus from the inequality \ref{eq15} we have:

\begin{equation}\label{eq4}
\frac{n}{4}\beta_X (\frac{R}{2}) \leq \beta_X ((m+2)R).
\end{equation}

We now use that $f(n)\leq C\beta_X(Cn+C)+C$ and  $\beta_X(n)\leq Cf(Cn+C)+C$, for all $n\geq n_0.$
From the first inequality and the fact that $\beta_X (Cn+C)\geq C$, we can deduce that 
\begin{equation}\label{eq5}
f(\frac{R}{4C}) \leq 2C \beta_X (\frac{R}{2}).
\end{equation}
 While from the second inequality, the fact that $R$ is big enough and the fact that $f(x)\geq kx$ (for $x\geq x_0$) we deduce that 

\begin{equation}\label{eq6}
\beta_X ((m+2)R)) \leq 2C f(4CmR).
\end{equation}

Combining inequalities \ref{eq4}, \ref{eq5} and \ref{eq6} we obtain

\begin{equation}\label{eq7}
 \frac{n}{4} \frac{1}{2C} f (\frac{R}{4C}) \leq 2C f(4CmR).
\end{equation}

Since $m>18$, from inequality \ref{eq13} we have that 

\begin{equation}\label{eq8}
 \frac{m}{2} R\leq  4 \phi
f^{-1}(n l_{n-2}),
\end{equation}

so $4CmR \leq  32 C \phi
f^{-1}(n l_{n-2})$, thus

\begin{equation}\label{eq9}
 f(4CmR) \leq  f(32 C) f(\phi) n l_{n-2}
\end{equation}

We then combine inequalities \ref{eq7} and \ref{eq9} and we obtain

\begin{equation}\label{eq10}
 \frac{n}{4} \frac{1}{2C} f (\frac{R}{4C}) \leq 2C f(32 C) f(\phi) n l_{n-2},
\end{equation}

then 

\begin{equation}\label{eq11}
 \frac{n}{4} \frac{1}{2C} k \frac{R}{4C} \leq 2C f(32 C) f(\phi) n l_{n-2},
\end{equation}

\begin{equation}\label{eq12}
 \frac{k}{64C^3} \frac{1}{f(32C)} R \leq  f(\phi)  l_{n-2}.
\end{equation}

Which is a contradiction since $f(\phi)< \frac{k}{64C^3} \frac{1}{f(32C)}$ and $l_{n-2} \leq R$.

\textbf{Case 2:} $\gamma \nsubseteq B_{d(x,y)}(x)$.\\
Then there exists $n_y \in \{ 1, \dots, n \} $ such that $x_j \in B_{d(x,y)}(x)$, for all $j \leq n_y-2$, while $x_{n_y -1} \notin B_{d(x,y)}(x)$. Then $d(x,y') \in (d(x,y)-R, d(x,y)]$, where $y'=x_{n_y -2}$.\\

Observe that the R-path $x=x_0, \dots,x_{n_y-2 }=y'$ belongs to $B_{d(x,y)} \cap A_i$, so we can use the arguments (slightly modified on constants) from Case 1. So the analogue of inequality \ref{eq13} is the following
 
 \begin{equation}\label{eq14}
(m-10)R\leq  4 \phi
f^{-1}(n l_{n-2})
\end{equation}

and since $m>20$, the inequality \ref{eq8} also holds for case 2.

 We continue as in case 1, and we reach the same contradiction.

\end{flushleft}

We set $\mathcal{U}_i$, to be the collection of the R-components of the annulus $A_i$, $i\geq2.$ Then the set $\mathcal{U}= B_{3Ri_0}(e) \cup \bigcup_{i\geq i_0}\mathcal{U}_i$, is a cover of $X$, with R-multipicity=2 ($X$ is unbounded and geodesic), and every element in $\mathcal{U}$ has diameter at most $\max \{6Ri_0, 30R \}$. Where $i_0=10[2R+4\phi f^{-1}(4R \beta_X((m+2)R))]$.

 It follows that $\operatorname{asdim} X = 1.$

\end{proof}

We can drop the requirement of $f$ being multiplicative, we have the following

\begin{theorem}
Let $X$ be an unbounded vertex-transitive connected graph. Suppose that there exist an invertible, non-decreasing,  $f$ and a constant $C>0$ such that $\beta_X(n) \simeq f(n)$ with constant $C$. We further suppose that $\beta_X(n) \succ n$. We assume that for every geodesic polygon with side lengths $l_0 \leq\cdots\leq l_{n}, $ we have 

$$
\operatorname{Slim}_X(P)
\leq \phi
f^{-1}(n  l_{n-2}),
$$

where $\phi \leq \frac{1}{32C}$. Then

$$
\operatorname{asdim}X=1.
$$

\end{theorem}
\begin{proof}
The proof of this theorem is almost identical to the proof of the previous. The only substantial difference is the following argument:\\ 
As previously, we have the following inequalities \begin{equation}\label{eq.14}
 \frac{n}{4} \frac{1}{2C} f (\frac{R}{4C}) \leq 2C f(4CmR),
\end{equation}

\begin{equation}\label{eq.15}
 4CmR \leq  32 C \phi
f^{-1}(n l_{n-2}).
\end{equation}

We then combine inequalities \ref{eq.14} and \ref{eq.15} and we obtain

\begin{equation}\label{eq.16}
 \frac{n}{8C}  f (\frac{R}{4C}) \leq 2C f( 32 C \phi
f^{-1}(n l_{n-2})  ) .
\end{equation}

We have that $\phi \leq \frac{1}{32C}$, and $l_{n-2}\leq R$, so

\begin{equation}\label{eq.17}
 \frac{n}{8C}  f (\frac{R}{4C}) \leq 2C f(
f^{-1}(n l_{n-2})  ),
\end{equation}
thus
\begin{equation}\label{eq.18}
  f (\frac{R}{4C}) \leq 16C  R .
\end{equation}

Which is a contradiction, since then we would have that $f(n)\preceq n$, by hypothesis we have that $f(n)\succ n$.

\end{proof}

As a corollary we have the following much cleaner result:
\begin{cor}
Let $X$ be an unbounded vertex-transitive connected graph. Suppose that there exist $r>1$ such that $\beta_X(n) \simeq n^r$ with constant $C$. We assume that for every geodesic polygon with side lengths $l_0 \leq\cdots\leq l_{n}, $ we have 

$$
\operatorname{Slim}_X(P)
\leq 
 \frac{1}{32C}\sqrt[r]{nl_{n-2} }.
$$
 Then

$$
\operatorname{asdim}_{}X=1.
$$
\end{cor}

\section{Trees of hyperbolic spaces}\label{section5}

\subsection{Trees of spaces}

\begin{defn}[Tree of spaces]\label{tree:def}
Let $T$ be a simplicial tree, with vertex set $V(T)$ and edge set $E(T)$.
A \emph{tree of spaces} over $T$ consists of the following data:
\begin{enumerate}
\item for every vertex $v\in V(T)$, a geodesic metric space $(X_v,d_v)$,
called the \emph{vertex space} of $v$;
\item for every edge $e\in E(T)$, a non-empty set $X_e$, called the
\emph{edge space} of $e$;
\item for every edge $e$ with endpoints $v$ and $w$, two maps
$$
f_e^v\colon X_e\to X_v, \qquad f_e^w\colon X_e\to X_w,
$$
called the \emph{attaching maps}.
\end{enumerate}
The \emph{total space} $X$ is obtained from the disjoint union
$$
\bigsqcup_{v\in V(T)}X_v\ \sqcup\ \bigsqcup_{e\in E(T)}X_e\times[0,1]
$$
by identifying $(x,0)$ with $f_e^v(x)$ and $(x,1)$ with $f_e^w(x)$, for
every edge $e$ with endpoints $v,w$ and every $x\in X_e$. For $x\in X_e$ the
segment $\{x\}\times[0,1]$ is called a \emph{rung} and has length $1$. We
equip $X$ with the induced length metric $d$: the distance between two
points is the infimum of the lengths of paths joining them which are
concatenations of geodesics of vertex spaces and of subsegments of rungs.
We always assume that $(X,d)$ is a geodesic metric space (this holds, for
instance, if all edge spaces are finite).

For an edge $e$ with endpoint $v$ we call
$$
A_e^v=f_e^v(X_e)\subseteq X_v
$$
the \emph{attaching set} of $e$ in $X_v$.
\end{defn}

\begin{rem}\label{tree:alt}
Some authors use a slightly different definition, in which the edge spaces
are subsets of the vertex spaces of their endpoints, $X_e\subseteq X_v\cap X_w$, 
and no rungs are added: the total space $X$ is obtained from the disjoint
union $\bigsqcup_{v\in V(T)}X_v$ by gluing vertex spaces along their edge spaces accordingly. This is the
case, for instance, for the Cayley graph of an amalgamated product over a
finite group. In this setting we additionally assume that distinct points
of attaching sets are at distance at least $1$ in every vertex space (this
holds, for instance, if the vertex spaces are graphs and the attaching sets
consist of vertices).
\end{rem}

\begin{defn}[Trees of hyperbolic spaces with bounded edge spaces]\label{tree:hyp}
Let $\delta\geq 0$ and $M\geq 0$. A tree of spaces $X$ is a \emph{tree of
$\delta$-hyperbolic spaces with $M$-bounded edge spaces} if
\begin{enumerate}
\item every vertex space $(X_v,d_v)$ is $\delta$-hyperbolic, and
\item every attaching set has diameter at most $M$ in its vertex space,
that is, $\operatorname{diam}_{d_v}(A_e^v)\leq M$ for every edge $e$ and
every endpoint $v$ of $e$.
\end{enumerate}
\end{defn}

\begin{exmp}\label{tree:exmp}
Let $G$ be the fundamental group of a finite graph of groups with finitely
generated vertex groups $G_v$ and finite edge groups $G_e$, and let $T$ be
its Bass--Serre tree. Taking as vertex space of the vertex $gG_v$ of $T$ the
coset $gG_v$ with the metric of the Cayley graph of $G_v$ (with respect to a
fixed finite generating set), as edge space of the edge $gG_e$ the finite
set $gG_e$, and as attaching maps the inclusions, we obtain a tree of
spaces with finite edge spaces. If all $G_v$ are hyperbolic, it is a tree
of $\delta$-hyperbolic spaces with $M$-bounded edge spaces, where $\delta$
and $M$ depend only on the finitely many vertex and edge groups.
\end{exmp}

In this section we prove the following theorem, which applies to both
definitions of trees of spaces.

\begin{theorem}\label{thm.trees of spaces}
Let $\delta\geq 0$ and $M\geq 0$. Assume that $X$ is a tree of
$\delta$-hyperbolic spaces with finite, $M$-bounded edge spaces. Assume
moreover that there are geodesic metric spaces $Y_1,\dots,Y_m$ and constants
$\lambda\geq1$, $\varepsilon\geq0$ such that every vertex space $X_v$ is
$(\lambda,\varepsilon)$-quasi-isometric to some $Y_j$, and every $Y_j$ is
$(\lambda,\varepsilon)$-quasi-isometric to some vertex space $X_v$. Then
$$
\operatorname{Slim}_X
\simeq
\max_{1\leq j\leq m}\operatorname{Slim}_{Y_j}.
$$
\end{theorem}

Note that each $Y_j$ is hyperbolic, being a geodesic metric space
quasi-isometric to a $\delta$-hyperbolic space.

\begin{rem}\label{tree:rem}
There are constants $\lambda'\geq1$, $\varepsilon'\geq0$ and $\delta'\geq0$,
depending only on $\delta$ and $M$, such that all the inclusion maps
$X_v\hookrightarrow X$ are $(\lambda',\varepsilon')$-quasi-isometric
embeddings, and $X$ is $\delta'$-hyperbolic. Indeed, a path in $X$ between
two points of $X_v$ which leaves $X_v$ through an attaching set $A_e^v$
must return to $X_v$ through the same attaching set; such an excursion has
length at least $1$ and can be replaced by a path in $X_v$ of length at most
$M$. This gives $d\leq d_v\leq\max\{1,M\}\,d$ on $X_v$. The hyperbolicity of
$X$ follows from the combination theorem of Bestvina and Feighn \cite{BF92},
since the edge spaces are uniformly bounded.
\end{rem}

\begin{flushleft}
\textbf{Note:} For the sake of readability, we give the proof in the case where $X$ is the
tree of spaces associated with a finite graph of groups with finite edge
groups (Example~\ref{tree:exmp}). The
proof of the general case is essentially the same, after replacing
cosets $gG_v$ of vertex groups by vertex spaces $X_v$, cosets $gG_e$ of edge
groups by attaching sets $A_e^v$, the Cayley graphs of the vertex groups by
the spaces $Y_j$, and the bound $\max_e\operatorname{diam}G_e$
by the constant $M$. The only properties used are that the vertex spaces
are uniformly hyperbolic and uniformly quasi-isometric to the spaces $Y_j$,
that the attaching sets are uniformly bounded, and Remark~\ref{tree:rem}.
\end{flushleft}



\subsection{Graphs of hyperbolic groups}
Let $(\mathbb{G}, Y)$ be a finite graph of groups with vertex groups $\lbrace G_{v} \mid v \in Y^{0} \rbrace$ and edge groups $\lbrace G_{e} \mid e \in Y^{1}_{+} \rbrace$.
Let $G=\pi_{1}(\mathbb{G},Y,\mathbb{T})$ be the fundamental group of this graph of groups (where $\mathbb{T}$ is a maximal subtree of $Y$).  We fix a finite generating set $S_v$ for the group $G_v$. Then the set $S=\cup S_v \bigcup \cup \lbrace  e , e^{-1} \mid e \notin E(\mathbb{T}))\rbrace$ is a finite generating set for the group $G$. We set $X_G=Cay(G,S)$ its Cayley graph.

The group $G$ acts on its Bass-Serre tree $T$. For every coset $g_vG_v$ there is
a natural projection $\pi_{G_{v}}: G \rightarrow T$ defined by the action: $ \pi_{G_{v}}(g)=  gg_vG_v$.

We consider the Bass-Serre tree $T$ as a metric space with the simplicial metric $d_T$. If $\Gamma$ is a graph, we denote by $\Gamma^{0}$, or $V(\Gamma)$ the vertices of $\Gamma$. 


\subsubsection{$\mathbb{G}$-decompositions}

A geodesic $n$-gon, say $P$, in $X_G$, is formed by its  vertices $V(P)=\{v_0, \dots, v_n\}$, and its geodesics $\gamma_i$ connecting $v_i$ to $v_{i+1}$. Let $p$ be a path in $X_G$, we denote by $p^{-1}$ the inverse path. We denote by $\partial P=\gamma$ the loop  $\gamma_0  \cdots  \gamma_n$.\\

\begin{defn}[VK-diagram]
Let $P$ be a geodesic n-gon, and $\mathbf{P}=\{P_i\}$ be a  finite family of geodesic polygons.  A \emph{VK-diagram} of $P$ over $\mathbf{P}$ is a finite, connected, simply connected planar $2$-complex $D$, together with an edge-labelling by elements of $S^{\pm1}$, satisfying the following conditions:

\begin{enumerate}
    \item Each oriented edge of $D$ is labelled by an element of $S^{\pm1}$. Reversing the orientation of an edge replaces its label by its inverse.

    \item The set of the words formed by the labels of the boundary path of each $2$-cell of $D$, read with a fixed orientation, is the set $\{\partial P_i\}$
\end{enumerate}

The label of the outer boundary path of $D$ is called the \emph{boundary label} of the diagram. We say that the collection $\mathbf{P}$ \emph{forms (is) a VK-diagram} of $P$.
\end{defn}
The idea is the following: we realise $P$ as the boundary of a 2-complex $D_i$, and  every $P_i$ as the boundary of a 2-cell $D$. Then $D$ is a VK-diagram of $P$ over $\mathbf{P}$ if the 2-cells of $D$ are exactly the $D_i$'s.\\

\begin{center}

\begin{tikzpicture}[scale=1.5]


\fill[blue!25]
  (0.9,0.15) .. controls (0.2,0.17) .. (-0.5,0.2)
  .. controls (-0.9,0.55) .. (-1.1,0.9)
  .. controls (-1.1,1.35) .. (-1.0,1.8)
  .. controls (-0.65,2.05) .. (-0.3,2.3)
  .. controls (0.1,2.2) .. (0.55,2.05)
  .. controls (0.57,1.5) .. (0.6,1.0)
  .. controls (0.75,0.57) .. (0.9,0.15)
  -- cycle;

\fill[red!25]
  (0.9,0.15) .. controls (0.75,0.57) .. (0.6,1.0)
  .. controls (0.57,1.5) .. (0.55,2.05)
  .. controls (0.85,2.3) .. (1.15,2.5)
  .. controls (1.5,2.0) .. (1.85,1.35)
  .. controls (1.95,1.25) .. (2.05,1.2)
  .. controls (1.47,0.67) .. (0.9,0.15)
  -- cycle;

\fill[green!30]
  (-0.3,2.3) .. controls (-0.15,2.7) .. (0.0,3.1)
  .. controls (0.35,3.35) .. (0.7,3.6)
  .. controls (1.05,3.57) .. (1.4,3.5)
  .. controls (1.55,3.25) .. (1.7,3.0)
  .. controls (1.3,2.75) .. (1.15,2.5)
  .. controls (0.85,2.3) .. (0.55,2.05)
  .. controls (0.1,2.2) .. (-0.3,2.3)
  -- cycle;

\fill[yellow!60]
  (2.05,1.2) .. controls (1.95,1.25) .. (1.85,1.35)
  .. controls (1.5,2.0) .. (1.15,2.5)
  .. controls (1.3,2.75) .. (1.7,3.0)
  .. controls (2.05,3.1) .. (2.4,3.2)
  .. controls (2.7,2.9) .. (3.0,2.6)
  .. controls (2.95,2.15) .. (2.9,1.7)
  .. controls (2.47,1.45) .. (2.05,1.2)
  -- cycle;

\draw[very thick]
  (0.9,0.15) .. controls (0.2,0.17) .. (-0.5,0.2)
  .. controls (-0.9,0.55) .. (-1.1,0.9)
  .. controls (-1.1,1.35) .. (-1.0,1.8)
  .. controls (-0.65,2.05) .. (-0.3,2.3)
  .. controls (-0.15,2.7) .. (0.0,3.1)
  .. controls (0.35,3.35) .. (0.7,3.6)
  .. controls (1.05,3.57) .. (1.4,3.5)
  .. controls (1.55,3.25) .. (1.7,3.0)
  .. controls (2.05,3.1) .. (2.4,3.2)
  .. controls (2.7,2.9) .. (3.0,2.6)
  .. controls (2.95,2.15) .. (2.9,1.7)
  .. controls (2.47,1.45) .. (2.05,1.2)
  .. controls (1.47,0.67) .. (0.9,0.15)
  -- cycle;

\draw[thick]
  (0.9,0.15) .. controls (0.75,0.57) .. (0.6,1.0)
  .. controls (0.57,1.5) .. (0.55,2.05);

\draw[thick]
  (-0.3,2.3) .. controls (0.1,2.2) .. (0.55,2.05)
  .. controls (0.85,2.3) .. (1.15,2.5);

\draw[thick]
  (1.15,2.5) .. controls (1.3,2.75) .. (1.7,3.0);

\draw[thick]
  (1.15,2.5) .. controls (1.5,2.0) .. (1.85,1.35)
  .. controls (1.95,1.25) .. (2.05,1.2);

\node at (-0.25, 1.1) {\large $P_1$};
\node at (0.95,  1.3) {\large $P_2$};
\node at (0.65,  2.9) {\large $P_3$};
\node at (2.2,   2.3) {\large $P_4$};
\node at (1.5,   0.0) {\large $P$};
\node at (3.5,   3.2) {\Large $D$};

\node[align=center] at (1.0,-0.5)
  {The 2-complex $D$ is a VK-diagram of $P$ over $\{P_1,\ldots,P_4\}$};

\end{tikzpicture}

\end{center}

\begin{defn}[$\mathbb{G}$-decompositions] Let $(\mathbb{G},Y)$ be a finite graph of groups, where $Y$ is a tree. A collection of finitely many geodesic polygons $\{ P_i \}$ is a \emph{$\mathbb{G}$-decomposition} of $P$ if the following five conditions are satisfied:
\begin{enumerate}
    \item the collection $\mathbf{P}=\{ P_i \}$ forms a VK-diagram of $P$.

    \item $V(P_i) \subseteq g_{u_{i}}G_{u_{i}}$
    \item every $P_i$ is an $n_i$-gon, where $n_i \leq n+2$
\item the vertices of $P_i$ are vertices of $P$, or elements contained in $gG_e \cap P$, for some edge groups.
\item every geodesic of $P_i$ is a geodesic segment of a geodesic of $P$, or a geodesic of $X_G$ connecting two elements in $gG_e \cap P$, which are contained in different geodesics of $P$, for some edge groups.

\end{enumerate}

\end{defn}

\begin{center}

\begin{tikzpicture}[scale=0.8]


%

\draw[blue, thick]
  (-0.35,-4.5)
  .. controls (-0.35,-3.5) .. (-0.35,-3.0)
  .. controls (-0.4,-2.2) and (-0.8,-1.5) .. (-1.3,-0.9)
  .. controls (-2.2, 0.2) and (-3.0, 1.2) .. (-3.3, 2.0)
  .. controls (-3.0, 2.8) and (-2.0, 2.8) .. (-1.7, 2.0)
  .. controls (-1.5, 1.4) and (-0.9, 0.5) .. (-0.5, 0.1)
  .. controls (-0.2,-0.1) and ( 0.2,-0.1) .. ( 0.5, 0.1)
  .. controls ( 0.9, 0.5) and ( 1.5, 1.4) .. ( 1.7, 2.0)
  .. controls ( 2.0, 2.8) and ( 3.0, 2.8) .. ( 3.3, 2.0)
  .. controls ( 3.0, 1.2) and ( 2.2, 0.2) .. ( 1.3,-0.9)
  .. controls ( 0.8,-1.5) and ( 0.4,-2.2) .. ( 0.35,-3.0)
  .. controls ( 0.35,-3.5) .. ( 0.35,-4.5)
  .. controls ( 0.2,-4.7) and (-0.2,-4.7) .. (-0.35,-4.5);

\filldraw[fill=blue!20, draw=none]
  (-0.35,-3.0)
  .. controls (-0.35,-3.5) .. (-0.35,-4.5)
  .. controls (-0.2,-4.7) and (0.2,-4.7) .. (0.35,-4.5)
  .. controls (0.35,-3.5) .. (0.35,-3.0)
  -- cycle;
\draw[red, line width=3pt] (-0.35,-3.0) -- (0.35,-3.0);
\node[blue, font=\large\bfseries] at (0,-3.9) {$P_1$};
\draw[red, line width=2pt] (0,-3.8) ellipse (2.2 and 0.9);
\node[red] at (3.3,-3.8) {$g_{u_0}G_{u_0}$};

\filldraw[fill=blue!20, draw=none]
  (-1.3,-0.9)
  .. controls (-2.2,0.2) .. (-1.825,-0.253)
  -- (-0.759, 0.394)
  .. controls (-0.9,0.5) .. (-0.5, 0.1)
  .. controls (-0.2,-0.1) and ( 0.2,-0.1) .. ( 0.5, 0.1)
  .. controls ( 0.9,0.5) .. ( 0.759, 0.394)
  -- ( 1.825,-0.253)
  .. controls ( 2.2,0.2) .. ( 1.3,-0.9)
  -- (-1.3,-0.9)
  -- cycle;

\draw[red, line width=3pt] (-1.825,-0.253) -- (-0.759, 0.394);
\draw[red, line width=3pt] ( 0.759, 0.394) -- ( 1.825,-0.253);
\draw[red, line width=3pt] (-1.3,-0.9) -- ( 1.3,-0.9);   

\node[blue, font=\large\bfseries] at (0,-0.3) {$P_2$};

\draw[red, line width=2pt] (0,-0.4) ellipse (2.5 and 0.8);
\node[red] at (3.6,-0.4) {$g_{u_1}G_{u_1}$};

\filldraw[fill=blue!20, draw=none]
  (-3.3, 2.0)
  .. controls (-3.0, 2.8) and (-2.0, 2.8) .. (-1.7, 2.0)
  -- (-3.3, 2.0) -- cycle;
\draw[red, line width=3pt] (-3.3,2.0) -- (-1.7,2.0);
\node[blue, font=\large\bfseries] at (-2.5, 2.5) {$P_4$};
\draw[red, line width=2pt] (-2.5, 2.3) ellipse (1.7 and 0.6);
\node[red] at (-2.5, 3.1) {$g_{u_2}G_{u_2}$};

\filldraw[fill=blue!20, draw=none]
  (1.7, 2.0)
  .. controls (2.0, 2.8) and (3.0, 2.8) .. (3.3, 2.0)
  -- (1.7, 2.0) -- cycle;
\draw[red, line width=3pt] (1.7,2.0) -- (3.3,2.0);
\node[blue, font=\large\bfseries] at (2.5, 2.5) {$P_3$};
\draw[red, line width=2pt] (2.5, 2.3) ellipse (1.7 and 0.6);
\node[red] at (2.5, 3.1) {$g_{u_3}G_{u_3}$};

\draw[blue, thick]
  (-0.35,-4.5)
  .. controls (-0.35,-3.5) .. (-0.35,-3.0)
  .. controls (-0.4,-2.2) and (-0.8,-1.5) .. (-1.3,-0.9)
  .. controls (-2.2, 0.2) and (-3.0, 1.2) .. (-3.3, 2.0)
  .. controls (-3.0, 2.8) and (-2.0, 2.8) .. (-1.7, 2.0)
  .. controls (-1.5, 1.4) and (-0.9, 0.5) .. (-0.5, 0.1)
  .. controls (-0.2,-0.1) and ( 0.2,-0.1) .. ( 0.5, 0.1)
  .. controls ( 0.9, 0.5) and ( 1.5, 1.4) .. ( 1.7, 2.0)
  .. controls ( 2.0, 2.8) and ( 3.0, 2.8) .. ( 3.3, 2.0)
  .. controls ( 3.0, 1.2) and ( 2.2, 0.2) .. ( 1.3,-0.9)
  .. controls ( 0.8,-1.5) and ( 0.4,-2.2) .. ( 0.35,-3.0)
  .. controls ( 0.35,-3.5) .. ( 0.35,-4.5)
  .. controls ( 0.2,-4.7) and (-0.2,-4.7) .. (-0.35,-4.5);

\node[blue] at (-4.8,-1.0) {\Large $P$};

\node[align=center] at (0,-5.8)
  {$\{P_1,\ldots,P_4\}$ is a $G$-decomposition of $P$};

\end{tikzpicture}
\end{center}

\begin{defn}[$\mathbb{G}_{HNN}$-decompositions] Let $(\mathbb{G},Y)$ be a finite graph of groups, where $Y$ is a rose. A collection of finitely many geodesic polygons $\{ P_i \}$ is a \emph{$\mathbb{G}_{HNN}$-decomposition} of $P$ if the following five conditions are satisfied:
\begin{enumerate}
    \item the collection $\mathbf{P}=\{ P_i \}$ forms a VK-diagram of $P$.

    \item $V(P_i) \subseteq N_1(g_{u_{i}}G_{u_{i}})$
    \item every $P_i$ is an $n_i$-gon, where $n_i \leq n+2$
\item the vertices of $P_i$ are vertices of $P$, or elements contained in $(gG_e\cup gt_{e}^{\pm 1}G_e) \cap P$ (where $t_e$ is the stable letter for the edge $e$), for some edge groups.
\item every geodesic of $P_i$ is a geodesic segment of a geodesic of $P$, or a geodesic of $X_G$ connecting two elements in $(gG_e\cup gt_{e}^{\pm 1}G_e) \cap P$, which are contained in different geodesics of $P$, for some edge groups.

\end{enumerate}

\end{defn}

We prove the following very important
\begin{lem}\label{lem5.1}
Let $P$ be n-gon in $X_G$,

\begin{itemize}
\item if $Y$ is a tree, then $P$ has a $\mathbb{G}$-decomposition,

\item if $Y$ is a rose, then $P$ has a $\mathbb{G}_{HNN}$-decomposition.
\end{itemize}
 
\end{lem}
\begin{proof} We prove only the first part, since the proof of the second part is essentially the same as the proof of the first part. Let $V(P)=\{v_0, \dots, v_n\}$ be the vertices of $P$, and let $\partial P=\gamma$ be the loop  $\gamma_0  \cdots  \gamma_n$. We set

\begin{equation}
\begin{aligned}
k(P)=\min\Bigl\{\, k_n+1 \ \Bigm|\;&
V(P)=\bigsqcup_{i=0}^{k_n}L_i,\quad
v_0\in L_0,\\
&L_0=\{v_{i_{-1}},\ldots,v_{i_0}\},\\
&\hspace{1.5cm}\ldots\\
&L_{k_n}=\{v_{i_{k_{n-1}}+1},\ldots,v_{i_{k_n}}\},\\
&i_{k_n+1}=i_{-1},\\
&\forall i\ \exists\, g_{u_i}G_{u_i}:
L_i\subseteq g_{u_i}G_{u_i},\\
&L_{i\pm1}\cap g_{u_i}G_{u_i}=\emptyset
\,\Bigr\}.
\end{aligned}
\end{equation}

This number roughly counts how many times the vertex set of \(P\) enters cosets of vertex groups of \(G\), taking multiplicity into account. The sets $L_i$ are disjoint. We will sometimes denote $g_{u_{i}}G_{u_{i}}$ by  $g_jG_j$.

We will use \emph{induction} on $k(P)$. For $k(P)=1$, there exists a coset of a vertex group, say $g_0G_0$, containing all the vertices of $P$. Then the set $\{P\}$ is a $\mathbb{G}$-decomposition of $P$.\\

For $k(P)=2$, then $V(P)=L_0 \sqcup L_1$. We will use another induction, on $d_T(v_{i_0 },v_{i_0 +1})= d_T(g_0G_0, g_1G_1)=\lambda_T(P)$.

For $\lambda_T=1$, we set $e=[g_0G_0,g_1G_1]$. Since $P$ enters $g_1G_1$, via $\gamma_{i_0 }$, it must exit $g_1G_1$, via $\gamma_{i_{1 } }$.
We recall that $\gamma_{i_0 }$ connects $v_{i_0 } \in L_0$ with $v_{i_{0} +1 } \in L_1$, and $\gamma_{i_{1 } }$ connects $v_{i_{1} }\in L_{i_{k_n } }$ with $v_{0} \in L_0$.
Obviously both those geodesics intersect a coset $g_eG_e$. Let $v' \in \gamma_{i_{1} } \cup g_eG_e$, and $v'' \in \gamma_{i_{0 } } \cup g_eG_e$. We consider a geodesic $r_e$ connecting $v'$ with $v''$.\\
Observe that the set of vertices $V_0=\{v',v''\}\cup L_0$, and $V_1=\{v',v''\}\cup L_1$, form geodesic $ n_i$-gons $P_0$ and $P_1$, respectively. Where $n_1, n_2 \leq n+1$. Observe that $k(P_0), k(P_1)=1$, since  $V(P_0)\subseteq g_0G_0$, and $V(P_1)\subseteq g_1G_1$. It is easy to see that $\{P_0, P_1\}$ is a $\mathbb{G}-$decomposition of $P$.

We assume is is true for $\lambda_T\leq m$, and let $P$ such that $\lambda_T=m+1$.\\
Let $w_0G_0, w_1G_1, \dots, w_mG_m$ a geodesic path from $g_0G_0$ to $g_1G_1$. We set $e=[w_0G_0, w_1G_1]$. Obviously $g_{i_0 } $ enters $g_1G_1$ via a coset $g_eG_e$, and $\gamma_{i_{1 } }$ exits $g_1G_1$ via $g_eG_e$.
Let $v' \in \gamma_{i_{1} } \cup g_eG_e$, and $v'' \in \gamma_{i_{0 } } \cup g_eG_e$. We consider a geodesic $r_e$ connecting $v'$ with $v''$.\\
Observe that the set of vertices $V_0=\{v',v''\}\cup L_0$, and $V_1=\{v',v''\}\cup L_1$, form geodesic $ n_i$-gons $P_0$ and $P_1$, respectively. Where $n_1, n_2 \leq n+1$. Observe that $k(P_0)=1$, $k(P_1)\leq2$, and $\lambda_T(P_1)\leq m$. Thus, the desired result follows by induction.\\

 We assume that the statement of the lemma is true for all geodesic polygons with $k(P)\leq k$. Let $P$ be a geodesic n-gon with $k(p)=k+1$.
 
So we have $V(P)=\bigsqcup_{i=0}^{k_n}L_i$. We may assume that $L_0=\{v_0,\dots ,v_{i_0 }\}.$ 

We use another induction on $\min \{ d_T(v_{i_0 },v_{i_0 +1})\} =\lambda_T$. We continue almost identically as the proof for the case $k(P)=2$. For $\lambda_T=1$, we find geodesic $n_i$-polygons $P_0$ and $P_1$, where $n_1, n_2\leq n+1$ and $k(P_1)=1$, while $k(P_0) \leq k$. Hence $P_0$ , and $P_1$ have $\mathbb{G}-$decompositions $\mathbf{P_0}$, and $\mathbf{P_1}$, respectively. The union of those decompositions is a $\mathbb{G}-$decomposition of $P$.

We assume is is true for $\lambda_T\leq m$, and let $P$ such that $\lambda_T=m+1$.\\
Again, we continue almost identically as the proof for the case $k(P)=2$. We find geodesic $ n_i$-gons $P_0$ and $P_1$, respectively. Where $n_1, n_2 \leq n+1$, and such that $k(P_0)=1$, $k(P_1)\leq k+1$, and $\lambda_T(P_1)\leq m$.  The desired result follows by induction.\\

\begin{flushleft}
\textbf{Note:} The definitions of $\mathcal{G}-$decompositions can be adjusted to the notion of trees of spaces. Furthermore, the lemma \ref{lem5.1} can be proved for the type of  trees of spaces considered in theorem \ref{thm.trees of spaces}.
\end{flushleft}



\end{proof}

\subsubsection{Main Result}
\begin{theorem}
Let $G$ be the fundamental group of a finite graph of hyperbolic groups, where every edge group is finite. 
Then 

$$\operatorname{Slim}_G
\simeq
\max_{v\in V(\mathcal{G})}\operatorname{Slim}_{G_v}.$$

\end{theorem}
\begin{proof}
Since all the egde groups are finite, we have that the vertex groups are quasi-convex subgroups of $G$. It follows that there exist constants $\lambda,\varepsilon >1$ such that the inclusion maps $i_v: G_v \hookrightarrow G$ are $(\lambda, \varepsilon)$-quasi-isometric embeddings, for all vertex groups. We conclude by Proposition \ref{prop3.1} that $\max \{ \operatorname{Slim}_{G_v}\} \preceq \operatorname{Slim}_G $.

We distinguish two cases.\\ 
\textbf{Case 1.} The graph $Y$ is a tree.\\
Let $P$ be a geodesic n-gon of $X_G$, then by lemma \ref{lem5.1} $P$ has a $\mathbb{G}-$decomposition, $\mathbf{P}=\{P_i \}$. Every $P_i$ is an $n_i$-gon, where $n_i \leq n+2$.\\
We recall that the vertices of $P_i$ are vertices of $P$, or elements contained in $gG_e \cap P$, for some edge groups.
We also recall that every geodesic of $P_i$ is a geodesic segment of a geodesic of $P$, or a geodesic of $X_G$ connecting two elements in $gG_e \cap P$, which are contained in different geodesics of $P$, for some edge groups.\\
 We observe then that any geodesic of $P_i$ which is a geodesic segment of $P$ is at distance at most $\max\{ \operatorname{Slim}_{(N_M(G_v),d_G)}(n+2)\}+ M$ from a different geodesic of $P$ =, where $M= \max\{\diam G_e\}$..

Using the fact that $\mathbf{P}=\{ P_i \}$ forms a VK-diagram of $P$, we have that any element of a geodesic of $P$ is contained in a geodesic of polygon $P_i$, for some $i$. We conclude that $\operatorname{Slim}_G(n) \leq \max\{ \operatorname{Slim}_{(N_M(G_v),d_G)}(n+2)\}+ M$.  So $$\operatorname{Slim}_G(n) \preceq \max\{ \operatorname{Slim}_{(N_M(G_v),d_G)}(n)\}$$
Since $(G_v,d_{G_v })$ and $(N_M(G_v),d_G)$, are quasi-isometric, using theorem \ref{thm3.1} we obtain the desired result.\\

\textbf{Case 2.} The graph $Y$ is not a tree.\\
Let $\mathbb{T}$ be a maximal tree of $Y$. 

We use induction on $n_e$, the number of edges of $Y\setminus T$. For $n_e=0$, it follows from case 1. We assume the the statement of the theorem is true for graphs with $n_e\leq m$, let $Y$ such that $n_e=m+1$.\\
Let $e$ be an edge of $Y\setminus T$. Let $Y_e$  be the graph obtained from  $Y$ by removing the edge $e$. Then observe that $G=\pi_1(\mathbb{G},Y,\mathbb{T})$ is an HNN-extension of $A=\pi_1(\mathbb{G},Y_e,\mathbb{T})$ over $G_e$.
Let $P$ be a geodesic n-gon of $X_G$, then by lemma \ref{lem5.1}, $P$ has a $\mathbb{G}_{HNN}-$decomposition, 
$\mathbf{P}=\{P_i \}$. 
Every $P_i$ is an $n_i$-gon, where $n_i \leq n+2$.

We recall that the vertices of $P_i$ are vertices of $P$, or elements contained in $(gG_e\cup gt_{e}^{\pm 1}G_e) \cap P$ (where $t_e$ is the stable letter for the edge $e$).
We further recall that, every geodesic of $P_i$ is a geodesic segment of a geodesic of $P$, or a geodesic of $X_G$ connecting two elements in $(gG_e\cup gt_{e}^{\pm 1}G_e) \cap P$, which are contained in different geodesics of $P$.\\
 We observe then, that any geodesic of $P_i$ which is a geodesic segment of $P$ is at distance at most $ \operatorname{Slim}_{(N_M(A),d_G)}(n+2)+ M+1$ from a different geodesic of $P$, where $M= \diam G_e$.

Using the fact that $\mathbf{P}=\{ P_i \}$ forms a VK-diagram of $P$, we have that any element of a geodesic of $P$ is contained in a geodesic of polygon $P_i$, for some $i$.  We conclude that $\operatorname{Slim}_{(N_M(A),d_G)}(n+2)+ \diam G_e+1$. So $$\operatorname{Slim}_G(n) \preceq \max\{ \operatorname{Slim}_{(N_M(A),d_G)}(n)\}.$$ 

Since $(A,d_{A })$ and $(N_M(A),d_G)$, are quasi-isometric, using theorem \ref{thm3.1} we obtain

$$\operatorname{Slim}_G(n) \preceq \max\{ \operatorname{Slim}_{(A,d_A)}(n)\}.$$ 

Then by induction we have we obtain the desired result.

\end{proof}

\section{Dichotomy for hyperbolic groups}
In this section we prove the following:

\begin{theorem}\label{dich:thmA}
Let $G$ be a finitely generated hyperbolic group.
\begin{enumerate}
\item If $G$ is virtually free, then $\operatorname{Slim}_G$ is bounded.
\item If $G$ is not virtually free, then $\operatorname{Slim}_G(n)\simeq\log n$.
\end{enumerate}
\end{theorem}

\subsection{A logarithmic upper bound}\label{dich:sec-upper}

\begin{prop}\label{dich:upper}
Let $X$ be a geodesic $\delta$-hyperbolic metric space. Then every
geodesic $n$-gon in $X$, $n\ge3$, is
$\delta\lceil\log_2(n-1)\rceil$-slim. In particular
$\operatorname{Slim}_X(n)\le\delta\lceil\log_2(n-1)\rceil$ and
$\operatorname{Slim}_X\preceq\log n$.
\end{prop}

The usual argument, triangulating the polygon by the diagonals from one
vertex, gives the bound $(n-2)\delta$. The improvement comes from
choosing a balanced triangulation, adapted to the side under
consideration. We first record an elementary identity which is used in
the induction step.

\begin{lem}\label{dich:ceil}
For every integer $s\ge2$,

$$
\big\lceil\log_2\lceil s/2\rceil\big\rceil=\lceil\log_2 s\rceil-1 ..
$$

In particular, if $s=2\lambda-1$ is odd, then
$$1+\lceil\log_2\lambda\rceil=\lceil\log_2 2\lambda\rceil=\lceil\log_2(2\lambda-1)\rceil.$$
\end{lem}

\begin{proof}
Recall that for integers $x\ge1$ and $k\ge1$ we have
$\lceil\log_2 x\rceil=k-1$ if and only if $2^{k-2}<x\le2^{k-1}$. Let
$k=\lceil\log_2 s\rceil$, so that $k\ge1$ and $2^{k-1}<s\le2^k$. It
suffices to show that $2^{k-2}<\lceil s/2\rceil\le2^{k-1}$.

If $s$ is even, then $\lceil s/2\rceil=s/2$ and the inequalities follow by
dividing by $2$.

If $s=2\lambda-1$ is odd, then $s\ge3$ and $\lceil s/2\rceil=\lambda$.
Since $\lambda>s/2>2^{k-2}$, the lower bound holds. For the upper bound,
note that $2^k$ is even while $s$ is odd, so $s\le2^k$ forces $s<2^k$,
that is $2\lambda=s+1\le2^k$. Hence $\lambda\le2^{k-1}$.

The last assertion follows since
$\lceil\log_2 2\lambda\rceil=1+\lceil\log_2\lambda\rceil$.
\end{proof}

\begin{proof}[Proof of Proposition~\ref{dich:upper}]
Let $P$ be a geodesic $n$-gon, and label its vertices $x_0,\dots, x_{n-1}$.
After relabelling cyclically, it suffices to show that the side
$[x_0,x_{n-1}]$ lies in the $\delta\lceil\log_2(n-1)\rceil$-neighbourhood
of the union of the other sides. For $1\le t\le n-1$ write
$e_t=[x_{t-1},x_t]$ for the sides of $P$, and for $0\le i<j\le n-1$ put

$$
A(i,j)=e_{i+1}\cup e_{i+2}\cup\dots\cup e_j .
$$

For every pair $i<j$ we fix a geodesic segment $\sigma_{ij}$ from $x_i$
to $x_j$ as follows: if $j=i+1$ we take $\sigma_{ij}=e_j$, we take
$\sigma_{0,n-1}=[x_0,x_{n-1}]$, the remaining side of $P$, for all other
pairs $\sigma_{ij}$ is an arbitrary geodesic segment.

\smallskip
\noindent\textbf{Claim.} For all $0\le i<j\le n-1$ and every
$q\in\sigma_{ij}$,

$$
d\big(q,A(i,j)\big)\le\delta\lceil\log_2(j-i)\rceil .
$$

We argue by induction on $s=j-i$. If $s=1$ then
$\sigma_{ij}=e_j\subseteq A(i,j)$ and there is nothing to prove. Let
$s\ge2$ and put $m=i+\lceil s/2\rceil$, so that $i<m<j$. The segments
$\sigma_{im}$, $\sigma_{mj}$, $\sigma_{ij}$ form a geodesic triangle,
which is $\delta$-slim. Hence there is
$q'\in\sigma_{im}\cup\sigma_{mj}$ with $d(q,q')\le\delta$. Both
$m-i=\lceil s/2\rceil$ and $j-m=\lfloor s/2\rfloor$ lie between $1$ and
$\lceil s/2\rceil$, and $A(i,m)\cup A(m,j)=A(i,j)$. By the induction
hypothesis,

$$
d\big(q',A(i,j)\big)\le\delta\big\lceil\log_2\lceil s/2\rceil\big\rceil=\delta\big(\lceil\log_2 s\rceil-1\big),
$$

where the last equality is Lemma~\ref{dich:ceil}. Therefore
$d(q,A(i,j))\le\delta+\delta(\lceil\log_2 s\rceil-1)=\delta\lceil\log_2 s\rceil$,
which proves the claim.


Applying the claim with $i=0$, $j=n-1$ gives the result, since
$A(0,n-1)$ is the union of all sides of $P$ other than $[x_0,x_{n-1}]$.
\end{proof}

\subsection{A logarithmic lower bound in the hyperbolic plane}\label{dich:sec-lower}

We will use two formulas for a geodesic right triangle in
$\mathbb{H}^2$ with legs $a,b$, hypotenuse $c$, and angle $\alpha$
opposite to $a$ (see e.g. \cite[\S3.5]{Ra}):
\begin{equation}\label{dich:trig}
\cosh c=\cosh a\cosh b,\qquad \sinh a=\sinh c\,\sin\alpha .
\end{equation}

\begin{prop}\label{dich:lower}
For every $n\ge3$ there is a geodesic $n$-gon $P_n$ in $\mathbb{H}^2$,
with a point $p$ on one of its sides, such that

$$
d\big(p,\ \text{union of the other sides of }P_n\big)\ \ge\ \log(n-1)-C_0,
\qquad C_0=\log\frac{\pi\sqrt5}{2}.
$$

Moreover, consecutive vertices of $P_n$ are at distance at most
$2\operatorname{arsinh}(1/2)<1$. In particular
$\operatorname{Slim}_{\mathbb{H}^2}(n)\ge\log(n-1)-C_0$.
\end{prop}

\begin{proof}
Fix $p\in\mathbb{H}^2$ and let $r_0>0$ be defined by
$\sinh r_0=(n-1)/\pi$. Let $\theta=\pi/(n-1)$ and let
$x_0,\dots,x_{n-1}$ be the points on the circle of radius $r_0$ about $p$
at angles $0,\theta,2\theta,\dots,(n-1)\theta=\pi$, measured from a fixed
direction at $p$. Then $x_0$ and $x_{n-1}$ are antipodal, so the geodesic
$[x_0,x_{n-1}]$ is a diameter and contains $p$. Let $P_n$ be the geodesic
$n$-gon with these vertices.

Consider a side $e=[x_{k-1},x_k]$, $1\le k\le n-1$, and let $\ell$ be its
length and $m$ its midpoint. The triangle with vertices $p,x_{k-1},x_k$ is
isosceles, so the geodesic from $p$ to $m$ is perpendicular to $e$ and
bisects the angle $\theta$ at $p$. Since the foot of the perpendicular
from $p$ to the geodesic line through $e$ is $m\in e$, the point of $e$
closest to $p$ is $m$. Put $r=d(p,m)$. Applying \eqref{dich:trig} to the
right triangle $p,m,x_k$ we get

$$
\sinh(\ell/2)=\sinh r_0\,\sin(\theta/2),\qquad \cosh r_0=\cosh r\,\cosh(\ell/2).
$$

Since $\sin(\theta/2)\le\theta/2=\pi/(2(n-1))$, the first identity gives
$\sinh(\ell/2)\le1/2$, hence $\ell\le2\operatorname{arsinh}(1/2)$ and
$\cosh(\ell/2)\le\sqrt5/2$. The second identity then gives

$$
e^{r}\ \ge\ \cosh r\ \ge\ \frac{2}{\sqrt5}\cosh r_0\ \ge\ \frac{e^{r_0}}{\sqrt5},
$$

so $r\ge r_0-\tfrac12\log5$. Finally
$r_0=\operatorname{arsinh}\big((n-1)/\pi\big)\ge\log\big(2(n-1)/\pi\big)$,
and therefore

$$
d(p,e)=r\ \ge\ \log(n-1)+\log\frac{2}{\pi}-\frac12\log5=\log(n-1)-C_0 .
$$

As this holds for every side $e\neq[x_0,x_{n-1}]$, the proposition
follows.
\end{proof}

\begin{prop}\label{dich:Hm}
For every $m\ge2$, $\operatorname{Slim}_{\mathbb{H}^m}(n)\simeq\log n$.
More precisely, for $n\ge3$,

$$
\log(n-1)-C_0\ \le\ \operatorname{Slim}_{\mathbb{H}^m}(n)\ \le\ \delta_0\lceil\log_2(n-1)\rceil,
\qquad \delta_0=\log(1+\sqrt2).
$$

\end{prop}

\begin{proof}
Geodesic triangles in $\mathbb{H}^m$ are $\delta_0$-slim, since every
triangle lies in a totally geodesic copy of $\mathbb{H}^2$, where the
optimal constant is $\log(1+\sqrt2)$. The upper bound follows from
Proposition~\ref{dich:upper}. For the lower bound, place the polygon $P_n$
of Proposition~\ref{dich:lower} in a totally geodesic copy of
$\mathbb{H}^2\subseteq\mathbb{H}^m$; the inclusion is isometric, so
distances are unchanged.
\end{proof}

\begin{rem}
For $\mathbb{H}^2$, Propositions~\ref{dich:upper} and~\ref{dich:lower}
give

$$
1\ \le\ \liminf_{n\to\infty}\frac{\operatorname{Slim}_{\mathbb{H}^2}(n)}{\log n}\ \le\
\limsup_{n\to\infty}\frac{\operatorname{Slim}_{\mathbb{H}^2}(n)}{\log n}\ \le\ \frac{\log(1+\sqrt2)}{\log2}\approx1.27 .
$$

\end{rem}

\subsection{Hyperbolic groups}\label{dich:sec-groups}
We are ready to prove Theorem~\ref{dich:thmA}
\begin{proof}[Proof of Theorem~\ref{dich:thmA}]
Fix a finite generating set of $G$ and let $\Gamma$ be the corresponding
Cayley graph, which is $\delta$-hyperbolic for some $\delta$.

(1) A finitely generated virtually free group is the fundamental group of
a finite graph of finite groups \cite{KPS}, so it acts properly and
cocompactly on its Bass--Serre tree $T$, which is locally finite. By the
\v{S}varc--Milnor lemma $\Gamma$ is quasi-isometric to $T$. Every geodesic
polygon in a tree is $0$-slim, so $\operatorname{Slim}_T\equiv0$, and by
Theorem~\ref{thm3.1} $\operatorname{Slim}_\Gamma\simeq\operatorname{Slim}_T$.
Hence $\operatorname{Slim}_\Gamma$ is bounded.

(2) The upper bound $\operatorname{Slim}_\Gamma\preceq\log n$ is
Proposition~\ref{dich:upper}. For the lower bound, since $G$ is not
virtually free, by \cite{BK} there is a
$(\lambda,\varepsilon)$-quasi-isometric embedding
$f\colon\mathbb{H}^2\to\Gamma$. By proposition \ref{prop3.1} and proposition \ref{dich:lower} we obtain the wanted result.

\end{proof}

\begin{proof}[Proof of the characterisation of virtually free groups]
(1)$\Rightarrow$(2) is Theorem~\ref{dich:thmA}(1), and (2)$\Rightarrow$(3)
is clear. Assume (3). Then $\operatorname{Slim}_\Gamma(3)<\infty$, so
geodesic triangles in $\Gamma$ are uniformly slim and $G$ is hyperbolic.
If $G$ is not virtually free, then by Bonk-Kleiner  \cite{BK} there is a quasi-isometric embedding
$f\colon\mathbb{H}^2\to\Gamma$.

Then by proposition \ref{prop3.1} and proposition \ref{dich:lower} we have
$\log n\preceq\operatorname{Slim}_\Gamma$, which contradicts (3)
\end{proof}


\begin{rem}
In view of \cite{PT Polygon} and proposition \ref{dich:upper}, if $G$ is a $\delta$ hyperbolic group, non-virtually free group,  with respect to a finite generating set with $k$ elements then we have the following very interesting inequality

\begin{equation}
  \frac{1}{8\log_{2}2k}\log_{2}(n\delta) <  \operatorname{Slim}_X(P) \leq 2 \log_2(n^\delta)
\end{equation}

for some $P$ n-gon (for infinitely many $n$) in the Cayley graph, $X$, of $G$. Indeed, if the left inequality failed for all but finitely many $n$, then,
since $\operatorname{Slim}_X(n)<\infty$ for every $n$, after enlarging
$\delta$ every geodesic $n$-gon in $X$ would be
$\frac18\log_{2k}(n\delta')$-slim for some $\delta'>0$, and $X$ would be a
quasi-tree by \cite{PT Polygon}, a contradiction.

Equivalently we have the following,
\begin{equation}\label{dich:ineq2}
n\delta\ <\ (2k)^{8\operatorname{Slim}_X(P)}\ \leq\ n^{16\delta\log_2(2k)} .
\end{equation}
Thus $(2k)^{8\operatorname{Slim}_X(P)}$ grows at least linearly and at most
polynomially in $n$.

\end{rem}

\section{Questions}\label{section 6}


\begin{flushleft}
\textbf{Question 2:} Can the assumption that the edge groups are finite in our graph-of-groups theorem be weakened?
More precisely, under what additional hypotheses on a finite graph of hyperbolic groups with possibly infinite edge groups can one still obtain
\[
\operatorname{Slim}_G
\simeq
\max_{v\in V(\mathcal G)}
\operatorname{Slim}_{G_v},
\]
where $G$ is the fundamental group of the graph of groups?
\end{flushleft}

\begin{flushleft}
\textbf{Question 3:} To what extent can the vertex-transitivity assumption in our results relating the growth of the slimness function to the asymptotic dimension be weakened or removed?
In particular, are there natural geometric conditions on a hyperbolic graph, that suffice to obtain analogous conclusions?
\end{flushleft}

The proof of Theorem~\ref{dich:thmA} uses the group structure only
through the theorem of Bonk and Kleiner. This suggests the following
questions.

\begin{flushleft}
\textbf{Question 4:}\label{dich:q1}
Let $X$ be a hyperbolic graph of bounded valency which is not a
quasi-tree. Is it true that $\operatorname{Slim}_X(n)\simeq\log n$? Does
this hold at least when $X$ is vertex-transitive?
\end{flushleft}

By \cite{PT Polygon} the answer to the second part is at least partially
positive: for vertex-transitive $X$ one has
$\operatorname{Slim}_X\not\preceq\sqrt{\log_k(n)}$.

\begin{flushleft}
\textbf{Question 5:}\label{dich:q2}
Is there a hyperbolic geodesic space, or a hyperbolic graph of bounded
valency, $X$ such that $1\preceq\operatorname{Slim}_X\preceq\log n$ with
$\operatorname{Slim}_X\not\simeq1$ and $\operatorname{Slim}_X\not\simeq\log n$?
\end{flushleft}

\begin{flushleft}
\textbf{Question 6:}\label{dich:q3}
Determine $\lim_{n\to\infty}\operatorname{Slim}_{\mathbb{H}^2}(n)/\log n$,
if it exists.
\end{flushleft}

\begin{flushleft}
Department of Mathematics University of the Aegean Karlovassi, 83200 Samos, Greece.\\
\textit{E-mail address}: ptselekidis@aegean.gr
\end{flushleft}


\begin{thebibliography}{99}






\bibitem{BK}
M. Bonk and B. Kleiner, {\em Quasi-hyperbolic planes in hyperbolic groups},
Proc. Amer. Math. Soc. 133 (2005), no. 9, 2491--2494.

\bibitem{BH99}
M. R. Bridson and A. Haefliger,
\textit{Metric Spaces of Non-Positive Curvature},
Grundlehren der mathematischen Wissenschaften, vol. 319,
Springer-Verlag, 1999.

\bibitem{GdlH}
\'E. Ghys and P. de la Harpe (eds.), {\em Sur les groupes hyperboliques
d'apr\`es Mikhael Gromov}, Progress in Mathematics, vol. 83,
Birkh\"auser, Boston, 1990.

\bibitem{Gr87}
M. Gromov, {\em Hyperbolic groups}, in Essays in Group Theory, Math. Sci.
Res. Inst. Publ., vol. 8, Springer, New York, 1987, 75--263.

\bibitem{KPS}
A. Karrass, A. Pietrowski, and D. Solitar, {\em Finite and infinite
cyclic extensions of free groups}, J. Austral. Math. Soc. 16 (1973),
458--466.

\bibitem{BF92} M. Bestvina and M. Feighn, {\em A combination theorem for negatively curved groups}, J. Differential Geom. 35 (1992), 85--101.


\bibitem{Ol91}
A. Yu. Ol'shanskii, {\em Periodic quotient groups of hyperbolic groups},
Mat. Sb. 182 (1991), no. 4, 543--567; English transl. in Math. USSR-Sb.
72 (1992), 519--541.

\bibitem{Ol93}
A. Yu. Ol'shanskii, {\em On residualing homomorphisms and $G$-subgroups
of hyperbolic groups}, Internat. J. Algebra Comput. 3 (1993), no. 4,
365--409.

\bibitem{Ra}
J. G. Ratcliffe, {\em Foundations of Hyperbolic Manifolds}, 2nd ed.,
Graduate Texts in Mathematics, vol. 149, Springer, New York, 2006.




\bibitem{HHS} Jason Behrstock, Mark F. Hagen, and Alessandro Sisto. Asymptotic dimension and smallcancellation for hierarchically hyperbolic spaces and groups. Proc. Lond. Math. Soc. (3),
114(5):890–926, 2017



\bibitem{BD04} G.Bell, and A.Dranishnikov, {\em On asymptotic dimension of groups acting on trees}.  Geometriae Dedicata, 103 (2004), 89-101





\bibitem{BD08} G. Bell and A. Dranishnikov. {\em Asymptotic dimension}. Topology
Appl., 155(12):1265-1296, 2008.

\bibitem{BF} Gregory C. Bell, Koji Fujiwara, {\em The asymptotic dimension of a curve graph is finite} , Journal of the London Mathematical Society, Volume 77, Issue 1, February 2008, Pages 33–50, https://doi.org/10.1112/jlms/jdm090
  
\bibitem{BDK} G. C. Bell, A. N. Dranishnikov, and J. E. Keesling, {\em On a formula for the asymptotic dimension of free products}, Fundamenta Mathematicae 183.1 (2004): 39-45.   
 
  

















 









\bibitem{Gr} M.Gromov, {\em Asymptotic invariants of infinite groups}, in Geometric Group Theory, v.
2, 1-295, London Math. Soc. Lecture Note Ser., vol. 182, Cambridge Univ. Press, Cambridge,
1993.







\bibitem{LynSch} R. Lyndon and P. Schupp, {\em Combinatorial group theory}, Springer-Verlag,
1977.











  

\bibitem{PT} Panagiotis Tselekidis, {\em Asymptotic Dimension of Graphs of Groups and One Relator Groups}. Algebraic and Geometric Topology 23 (2023) 3587–3613. DOI: 10.2140/agt.2023.23.3587





\bibitem{PT Polygon} Panagiotis Tselekidis, {\em Logarithmically Slim Polygons}, (2026), Arxiv: 2611.05429





\bibitem{Yu} Guoliang Yu, {\em The Novikov conjecture for groups with finite asymptotic dimension},



\end{thebibliography}
\end{document}